# Finding New Limit Points of Mahler Measure by Methods of Missing Data Restoration

**Jean-Marc Sac-Épée**
Research engineer, IECL, University of Lorraine, France
jean-marc.sac-epee@univ-lorraine.fr

**Souad El Otmani**
Associate Professor, IUT de Metz, France
souad.elotmani@univ-lorraine.fr

**Armand Maul**
Professor, LIEC, University of Lorraine, France
armand.maul@univ-lorraine.fr

**Georges Rhin**
Professor, IECL, University of Lorraine, France
georges.rhin@univ-lorraine.fr

## Abstract

It is well known that the set of Mahler measures of single variable polynomial has limit points of which a list established by D. Boyd and M. Mossinghoff has been extended through approaches based on genetic algorithms. In this paper, we wish to further extend the list of known limit points by adapting a method of missing data restoration.

Keywords: Polynomials, Mahler measure, restoration methods, algorithms, number theory

# 1. INTRODUCTION

In 1962 K. Mahler defined (Mahler, 1962) what was later called the Mahler measure of a polynomial as follows: let $P$ be a polynomial with complex coefficients. $P$ can be written as

$$P(X)=a_0X^d+a_1X^{d-1}+\cdots+a_d=a_0\prod_{i=1}^{d}(X-\alpha_i),$$

where $a_0 \neq 0$ is the leading coefficient of P and the $\alpha_i$'s are its complex roots. Then, the Mahler measure of $P$ is defined as

$$M(P):=|a_0|\prod_{i=1}^{d}\max(1,|\alpha_i|).$$

In the following, we summarize some important properties of the Mahler measure. For a complete overview of the classical and more recent results concerning the Mahler measure, we refer to (Smyth, 2008). Recall that the Mahler measure $M(\alpha)$ of an algebraic number $\alpha$ is simply defined as that of its minimal polynomial $P_\alpha$ in $Z[X]$.

**Proposition 1.** If $\alpha$ is an algebraic number, then $M(\alpha) \geq 1$. Moreover, $M(\alpha) = 1$ if and only if $\alpha$ is a root of unity.

An algebraic number is said to be *reciprocal* if its minimal polynomial is reciprocal, i.e. if $P_\alpha(X)=X^dP_\alpha(1/X)$.

In (Smyth, 1971), C. Smyth proves the following result:

**Proposition 2.** Let $\theta_0$ be the only real root of the equation $\theta^3-\theta-1=0$. If $\beta$ is an algebraic integer such that $M(\beta)<\theta_0$, then $\beta$ is reciprocal.

From Proposition 1, we know that an algebraic number is a root of unity if and only if its Mahler measure is 1. In (Lehmer, 1993), Lehmer then wonders whether one can approach 1 as closely as desired by Mahler measures of algebraic integers which are not roots of unity. This problem is, to date, an open problem. Lehmer gave the smallest known Mahler measure > 1, $M(P_0) = 1.176280...$, where

$$P_0(X)=X^{10}+X^9-X^7-X^6-X^5-X^4-X^3+X+1.$$

M. Mossinghoff collected on a website (Mossinghoff) all known polynomials (monic, irreducible in $Z[X]$, with integer coefficients) of Mahler measure below 1.3, and the corresponding measures. Our purpose here is to deal with the problem of finding new limit points of Mahler measures. For this, let us first recall that the Mahler measure of a multivariable polynomial $P$ (in $n$ variables) is defined as

$$M(P)=\exp\left(\int_0^1\ldots\int_0^1\log\left|P\left(\exp(2i\pi t_1),\ldots,\exp(2i\pi t_n)\right)\right|dt_1\ldots dt_n\right).$$

In one dimension, the previous definition $M(P):=|a_0|\prod_{i=1}^{d}\max(1,|\alpha_i|)$ is derived from this one by using Jensen's formula.

Now, the point is that measures of multivariable polynomial are limiting values of measures of polynomials in fewer variables (Boyd et al., 2005, Lawton, 1983, Boyd, 1981). In (Boyd et al., 2005), the authors give 48 irreducible polynomials in two variables with Mahler measures below 1.37. The Mahler measures of these polynomials are limit of Mahler measures of univariate polynomials. In (El Otmani et al., 2017), the authors add 11 new polynomials to the previous list using an approach based on genetic algorithms.

In this work, our aim is to extend these lists by implementing a method based on the technics of missing data restoration, which are frequently used in different fields of mathematics, for example in tomographic reconstruction.

## 2. MISSING DATA RESTORATION ALGORITHM

### 2.1 Motivation of The Approach

The idea of this attempt is to try to discover new bivariate polynomials with Mahler measure below 1.37, not directly as in (Boyd et al., 2005) and (El Otmani et al., 2017), but by trying to exploit lists of already known polynomials. To summarize the approach in a few words, the principle of the method consists in assigning to a polynomial a small number of random coefficients, and in considering its other coefficients as missing data. These missing coefficients are then reconstructed by a restoration algorithm, which replaces them with the coefficients of the closest known polynomial in the sense of a distance to be defined. For the method to have a chance of success, it is necessary to have a sufficient number of polynomials of the same type, i.e., in particular, with the same degree. So we started by taking the 11 polynomials given in (El Otmani et al., 2017), because these polynomials are all of the same type, namely polynomials of the form $P_{a_0,\ldots,a_{35}}(x,y)=a_0+a_1y+a_2y^2+a_3y^3+a_4y^4+a_5y^5+a_6x+a_7xy+a_8xy^2+a_9xy^3+\cdots+a_{34}x^5y^4+a_{35}x^5y^5$, with $a_0,a_1,\ldots,a_{35}$ in $\{-1,0,1\}$. It can be noted that for each of the 11 polynomials in the list, only 6 of its 36 coefficients are non-zero. This is not particularly important for our approach. We then added to this list those polynomials provided in (Boyd, 2005) with degrees low enough to be written in the form specified above simply by adding zeros for some coefficients. We have thus obtained a list of reference polynomials for use in our algorithm.

In the following lines, we outline the principle of restoring missing data.

### 2.2 The Principle of The Restoration Algorithm

In our case, the principle of the restoration algorithm is to consider the polynomials of the reference list as vectors with 36 coefficients. To evaluate how these vectors, relate to each other, a distance is introduced. We have a wide choice of distances available (Mahalanobis distance, Minkowski distance, Hamming distance, Chebychev distance...) among which we choose the most usual one, namely the euclidean distance $\mathrm{dist}(P_1,P_2)=\sqrt{\sum_{i=1}^{36}\left(x_i^2-x_i^1\right)^2}$, where $P_k=\left(x_1^k,\ldots,x_{36}^k\right)$, k=1,2. We can set the leading coefficient of the unknown polynomial to 1, so there are still $3^{35}=50031545098999707$ coefficients to find, which makes it impossible to try all the possibilities, or to try to find them by random draws.

Our algorithm works as follows:

1. We randomly draw an integer $N$ between 1 and 10.
2. $N$ coefficients are randomly selected among the 36 coefficients of a polynomial which is of the type specified in section 2.1.
3. Each of these $N$ coefficients are assigned a random value selected from $\{-1;0;1\}$.
4. The distances between this polynomial with missing data and each of the polynomials in the reference list are computed. Only the assigned coefficients are used to compute the distances.
5. The missing coefficients of the polynomial are replaced by the corresponding coefficients of the nearest polynomial (for the Euclidean distance).
6. We control the Mahler measure and the irreducibility of the obtained polynomial. If its Mahler measure (or that of one of its factors) is less than 1.37, it is compared to that of already known polynomials, and kept if it is not listed.
7. Back to 1.

## 3. IMPLEMENTATION AND RESULTS

For the concrete implementation of our restoration algorithm, we have used the Matlab Statistic Optimization Toolbox (Coleman) which even allows, in a broader way, to rely on reference vectors that themselves have missing data. In the case we deal with, of course, we only use a list of fully known vectors. As in (El Otmani et al., 2017), the irreducibility of the polynomials in two variables was verified using the symbolic-based mathematical software Maxima, and the accurate calculation of the measures of two-variable polynomials was performed using GNU Octave.

The calculations, which have lasted for 2 weeks, were performed on a Dell Precision M6700 (processor: Intel Core i7-3940XM CPU @ 3.00 Ghz x 8, memory: 15.6 Gio).

Note that our algorithm gave us limit points already provided in (Boyd et al., 2005), whose corresponding polynomials did not appear among our reference polynomials because their degrees did not allow them to be written in the form

$$P_{a_0,\dots,a_{35}}(x,y)=a_0+a_1y+a_2y^2+a_3y^3+a_4y^4+a_5y^5+a_6x+a_7xy+a_8xy^2+a_9xy^3+\cdots+a_{34}x^5y^4+a_{35}x^5y^5.$$

For example, the point 1.358545590 of (Boyd et al., 2005) which corresponds to polynomial $P(4,1)=x^6+x^5+x^4+x^3y^2+x^3y+x^3+x^2y^2+xy^2+y^2$ is found by our algorithm via polynomial $x^2y^4-x^2-xy^3-xy^2-y^5+y$. The point 1.366145966 of (Boyd et al., 2005) which corresponds to polynomial

$P(5,3)=x^6+x^5+x^4y^2+x^4y+x^4+x^3y^2+x^3y+x^3+x^2y^2+x^2y+x^2+xy^2+y^2$ is found by our algorithm via polynomial $x^2y^5+x^2+xy^4+xy+y^5+1$. The point 1.366807889 of (Boyd et al., 2005) which corresponds to polynomial $P(5,1)=x^8+x^7+x^6+x^5+x^4y^2+x^4y+x^4+x^3y^2+x^2y^2+xy^2+y^2$ is found by our algorithm via polynomial $x^2y^5-x^2-xy^3+xy^2+y^5-1$. Of course, these polynomials were immediately added to our list of reference polynomials to help find still unknown limit points.

The table below gives the new limit points provided by our algorithm.

| **Mahler measures** | **Polynomials** |
|---|---|
| 1.359375641 | $1-y-x^2-x^3y^3-x^5y^2+x^5y^3$ |
| 1.368922213 | $y-y^3-x^2-x^3y^4-x^5y+x^5y^3$ |

## 4. FINAL REMARKS

Data recovery algorithms may be criticized for not being able to produce new results without starting from a sufficient number of already known results. On the other hand, they offer the possibility to quickly extend an already existing list by exploiting the polynomials of this list, which is an alternative approach that has proved effective in practice.